\documentclass[11pt]{article}
\usepackage{amsfonts}
\usepackage{amsthm}
\usepackage{amsmath}
\usepackage{amssymb}
\usepackage{amscd}
\usepackage{verbatim}
\usepackage{multicol}

\theoremstyle{definition}
\newtheorem{thm}{Theorem}[section]
\newtheorem{prop}[thm]{Proposition}

\newtheorem{lem}[thm]{Lemma}

\newtheorem{ex}{Example}[section]

\numberwithin{equation}{section}

\def\obs{\noindent{\normalsize{\bf Remark: }}}

\def\blb#1{\text{$\mathbb{#1}$}}
\def\cal#1{\text{$\mathcal{#1}$}}

\def\parcial#1/#2^#3{\frac{\partial^{#3}#1}{\partial#2^{#3}}}
\def\ord#1^#2{#1$^{\text{#2}}$}
\def\foral{\,\forall\,}
\def\parcial#1/#2^#3{\frac{\partial^{#3}#1}{\partial#2^{#3}}}
\def\parcialc#1/#2+#3{\frac{\partial^2#1}{\partial#2 \partial#3}}
\def\parcialf#1/#2{\frac{\partial#1}{\partial#2}}
\def\lie#1{\mathfrak{#1}}

\def\hlie#1{\widehat{\mathfrak{#1}}}
\def\kz{Knizhnik-Zamolodchikov }
\def\repf{\text{\bf Rep}_{\text{f}}}

\begin{document}


\author{Pavel I. Etingof\hspace{1.35cm} Adriano A. Moura
\\ {\small etingof@math.mit.edu \hspace{.5cm} adrianoam@ime.unicamp.br }}
\title{On the quantum Kazhdan-Lusztig functor}
\date{March/2002}
\maketitle


\section*{Introduction} 

\nocite{ABRR98}
\nocite{EFK98}
\nocite{EtSc98b,EtSc99}
\nocite{EtVa98,EtVa99,EtVa00}
\nocite{Fel94a,Fel94b}
\nocite{FeVa96}
\nocite{Fra02}
\nocite{FrRe92}
\nocite{JKOS99}
\nocite{KaLu93a,KaLu93b,KaLu94a,KaLu94b}
\nocite{KaSo95}
\nocite{Mou01}
\nocite{TaVa97,TaVa01}

One of the most exciting developments in representation theory in the recent years
was the discovery of the Kazhdan-Lusztig functor \cite{KaLu93a,KaLu93b,KaLu94a,KaLu94b},
which is a tensor functor from the fusion category of representations 
of an affine Lie algebra to the category of representations of the corresponding 
quantum group, and is often an equivalence of categories. 
Informally speaking, this functor may be obtained by 
considering the monodromy of the Knizhnik-Zamolodchikov equations, or, 
equivalently, by considering exchange relations between 
intertwining (vertex) operators for the affine Lie algebra. 

Therefore the introduction of quantized vertex operators and quantized 
Knizhnik-Zamolodchikov equations by Frenkel and Reshetikhin \cite{FrRe92} gave rise to 
a hope that this remarkable functor may be q-deformed. In fact, it was 
implicitly suggested already in \cite{FrRe92} that, since monodromy matrices of 
the quantized KZ equations are elliptic, the q-deformed KL functor should 
map the category of representations of the quantum affine algebra
to the category of representations of an appropriate elliptic quantum group. 
However, at that time it was not clear what exactly 
this elliptic quantum group should be. This became much more clear after the works
\cite{Fel94a,Fel94b,FeVa96,TaVa97}. Namely,
in \cite{Fel94a,Fel94b} elliptic quantum groups were defined 
using a dynamical R-matrix, 
in \cite{FeVa96} their representations were studied
in the $sl(2)$ case, and in \cite{TaVa97} it was shown that 
the monodromy of the quantum KZ equations for $sl(2)$ 
is given by the R-matrices for these representations
(this is a q-analog of the Drinfeld-Kohno theorem). More specifically,
since these works and especially the later papers \cite{JKOS99},
\cite{ABRR98},\cite{EtVa99} it became transparent that the elliptic quantum group
may be obtained by twisting the usual quantum affine algebra by an appropriate 
dynamical (or quasi-Hopf) twist\footnote{We note that the idea of the
dynamical twisting first appeared in the important paper
\cite{BBB95}.}. 
In other words, there exists a tensor functor 
from the category of finite dimensional representations of the quantum affine algebra 
to the category of finite dimensional representations of the elliptic quantum group, 
which could be thought of as a q-analog of the Kazhdan-Lusztig
functor (See Section 5). This idea gave rise 
to many interesting works, where the appropriate twist 
and its limiting cases were carefully studied and computed 
(see e.g. \cite{Fra02} and references therein). Nevertheless, 
to the best of our knowledge, 
the corresponding quantum Kazhdan-Lusztig 
functor was never explicitly studied, and 
in particular it was never shown that it is fully faithful (i.e. defines 
an equivalence of the source category with some full subcategory
of the target category). 

The goal of this paper is to give a rigorous definition and to study this 
(essentially, previously known) functor, following 
the ideas of \cite{JKOS99},\cite{EtVa99}. Namely, our main result can 
be stated as follows.

\noindent{\bf Theorem.} There exists an exact and fully faithful 
tensor functor $\cal F$ from the category of finite dimensional 
representations of the quantum affine algebra $U_q({\hlie sl}_n)$ 
(with extended scalars)
to the category of finite dimensional representations of Felder's
elliptic quantum group $E_{\tau,\gamma/2}({\hlie gl}_n)$, for $q=e^{\pi i\gamma}$. 

This theorem implies that
 the category of representations of the quantum affine algebra can be viewed 
as a full subcategory of the category of representations of the elliptic quantum group. 
Furthermore, it can be shown 
that evaluation representations of the quantum affine algebra 
are mapped (up to tensoring with a 1-dimensional representation) 
to the corresponding 
evaluation representations of the elliptic quantum group, 
which are constructed in \cite{TaVa01} by pulling back 
finite dimensional representations of the 
small elliptic quantum group (with standard dynamical highest weight).
Thus, the image of $\cal F$ contains (up to tensoring with
1-dimensional representations)
all subquotients of tensor products 
of evaluation representations, and they have the same characters as the corresponding 
subquotients of tensor products of evaluation representations of the 
quantum affine algebra.  Moreover, it is easy to extend the functor $\cal F$ 
to the category of possibly infinite dimensional representations of 
$U_q({\hlie sl}_n)$, which belong to category 
$\cal O$ as representations of $U_q({\lie sl}_n)$. 
The existence of such extension implies that all such representations 
of $U_q({\hlie sl}_n)$ admit an elliptic deformation. This confirms 
a conjecture from \cite{TaVa01}, where these elliptic deformations 
were constructed for Verma and finite dimensional modules, and conjectured to exist for any 
highest weight module. 

\noindent{\bf Acknowledgments.}
P.E. is grateful to A. Varchenko for many useful discussions about 
the subject matter of this paper, and to G. Felder for explaining
the content of the work \cite{Ca01}. The work of P.E. was partially 
supported by the NSF grant DMS-9988796, and was done in part for the Clay Mathematics 
Institute. A.M. is grateful to MIT for hospitality. The 
Ph.D. studies of A.M., during which this work was accomplished,  
are supported by FAPESP (99/11600-0), and his visit to MIT
is supported by CAPES (0365/01-7), Brazil.


\section{Representations of Quantum Dynamical R-Matrices with Spectral Parameters}

Let $\lie a$ be a finite dimensional abelian Lie algebra over \blb C,
$V$ be a finite dimensional semisimple $\lie a$-module,
and $\gamma \in \blb C^*$. The quantum dynamical Yang-Baxter equation
with spectral parameter and step $\gamma$ is
\begin{align}\notag
R^{12}(u_1-u_2,\lambda-\gamma h^{(3)}) R^{13}(u_1-u_3,\lambda) R^{23}(u_2-u_3,\lambda-\gamma h^{(1)})\\ \label{eq:qdybs}
= R^{23}(u_2-u_3,\lambda) R^{13}(u_1-u_3,\lambda-\gamma h^{(2)}) R^{12}(u_1-u_2,\lambda)
\end{align}
with respect to a meromorphic function $R:\blb C\times\lie a^* \to \text{End}_{\lie a}(V\otimes V)$. The notation $h^{(i)}$ means that, for example, $R^{12}(u_1-u_2,\lambda-\gamma h^{(3)})(v_1\otimes v_2\otimes v_3) = \big(R(u_1-u_2,\lambda-\gamma\mu)(v_1\otimes v_2)\big)\otimes v_3$ if $v_3$ has weight $\mu$. A function $R$ is called a quantum dynamical R-matrix if it is a generically invertible solution of \eqref{eq:qdybs}.

Given a quantum dynamical R-matrix $R:\blb C\times\lie a^* \to \text{End}_{\lie a}(V\otimes V)$, we denote by $\repf(R)$ the category of meromorphic finite dimensional representations of $R$. This category was first defined in \cite{Fel94a,Fel94b}. It consists of pairs $(W,L_{W})$ where $W$ is a semisimple finite dimensional $\lie a$-module and $L_{W}$ is a meromorphic function
$$L_{W} : \blb C \times \lie a^* \to \text{End}_{\lie a}(V\otimes W)$$
such that $L_{W}(u,\lambda)$ is (generically) invertible and satisfies
\begin{align}\notag
R^{12}(u_1-u_2,\lambda-\gamma h^{(3)}) L_{W}^{13}(u_1-u_3,\lambda) L_{W}^{23}(u_2-u_3,\lambda-\gamma h^{(1)})\\ \label{eq:Ldef}
= L_{W}^{23}(u_2-u_3,\lambda) L_{W}^{13}(u_1-u_3,\lambda-\gamma h^{(2)}) R^{12}(u_1-u_2,\lambda)
\end{align}                         
in $\text{End}_{\lie a}(V\otimes V\otimes W)$. If $(W,L_{W})$ and $(U,L_{U})$ are two representations, a morphism between them is a meromorphic function $f:\lie a^* \to \text{Hom}_{\lie a}(W,U)$ such that
\begin{equation}\label{eq:morpdef}
\big(1\otimes f(\lambda)\big)L_{W}(u,\lambda) = L_{U}(u,\lambda)\big(1\otimes f(\lambda - \gamma h^{(1)})\big)
\end{equation} 
Thus, the category $\repf(R)$ is linear over the field $M_V$ of 
meromorphic functions on $\lie a^*$ periodic under the lattice spanned by 
the weights of $V$. 

\begin{ex}\label{ex:bex} Let $R: \blb C\times \lie a^* \to \text{End}_{\lie a}(V\otimes V)$ be a quantum dynamical R-matrix.
\begin{enumerate}\vspace{-.25cm}
\item\label{ex:trep}
Let $W$ be any vector space and let $\lie a$ act trivially on $W$. Then $(W,1)$, where $1$ denotes the constant function $1\in \text{End}_{\lie a}(V\otimes W)$, is a representation  of $R$ called the trivial representation associated to $W$.  

\item $(V,R)$ is a representation of $R$ called the basic representation.

\item\label{ex:isorep}
If $(W,L_{W})$ is a representation of $R$, $f:\lie a^* \to \text{End}_{\lie a}(W)$ is generically invertible and
$L_{W}^f = \big(1\otimes f(\lambda)^{-1}\big)L_{W}(u,\lambda)\big(1\otimes f(\lambda -\gamma h^{(1)})\big)$
then  $(W,L_{W}^{f})$ is also a representation. It is clear that $f$ provides an isomorphism from $(W,L^f_W)$ to $(W,L_W)$.
\end{enumerate}  
\end{ex}

\begin{ex}\label{ex:felder}
Let $\lie a$ be an abelian $n$-dimensional Lie algebra and $V$ be a semisimple $n$-dimensional $\lie a$-module whose weights form a basis for $\lie a^*$
(in other words, $\lie a$ is the usual Cartan subalgebra of $\lie{gl}_n$ acting naturally on $V=\blb C^n$). Quantum dynamical R-matrices for such $\lie a$ and $V$ are said to be of  $\lie{gl}_n$ type \cite{EtVa98,EtVa99}. They have the form
$$R(u,\lambda) = \sum_{m,l}  \alpha_{m,l}(u,\lambda) E_{m,m}\otimes E_{l,l} + \sum_{m\neq l} \beta_{m,l}(u,\lambda) E_{l,m}\otimes E_{m,l}$$
The main example is 
 the elliptic quantum dynamical R-matrix found by Felder \cite{Fel94a,Fel94b} 
$$R_{\tau,\gamma}^{ell}(u,\lambda) = \sum_{m} E_{m,m}\otimes E_{m,m} + \sum_{m\neq l} \alpha(u,\lambda_{m,l}) E_{m,m}\otimes E_{l,l} + \beta(u,\lambda_{m,l}) E_{l,m}\otimes E_{m,l}$$
with
$$\alpha(u,\lambda) = \frac{\vartheta_1(\lambda+\gamma;\tau)}{\vartheta_1(\lambda;\tau)} \frac{\vartheta_1(u;\tau)}{\vartheta_1(u-\gamma;\tau)}\quad\quad \beta(u,\lambda) = \frac{\vartheta_1(\gamma;\tau)}{\vartheta_1(\lambda;\tau)} \frac{\vartheta_1(u-\lambda;\tau)}{\vartheta_1(u-\gamma;\tau)}$$
where $\lambda_{m,l} = \lambda_m - \lambda_l$, and $\vartheta_1$ is the standard first theta function. Felder also proposed to associate to $R^{ell}_{\tau,\gamma}$ 
an elliptic quantum group $E_{\tau,\gamma/2}({\lie gl}_n)$, 
whose category of finite dimensional representations is $\repf(R^{ell}_{\tau,\gamma})$.
\end{ex}

Let $R: \blb C\times \lie a^* \to \text{End}_{\lie a}(V\otimes V)$ be a quantum dynamical R-matrix. The category $\repf(R)$ can be endowed with a tensor product functor given by
$(W,L_{W})\odot (U,L_{U}) = (W\otimes U,L_{W\odot U})$ where
$$L_{W\odot U}(u,\lambda) = L_{W}^{12}(u,\lambda-\gamma h^{(3)})L_{U}^{13}(u,\lambda)$$
on objects and
$$f\odot g(\lambda) = f(\lambda-\gamma h^{(2)})\otimes g(\lambda)$$
on morphisms. One easily checks that the first formula defines a representation, 
and the second one a morphism of representations. 
Thus, $\odot$ is a bifunctor. Moreover, it is easy to see that 
it defines a structure of a tensor category on $\repf(R)$, 
with unit object being the trivial representation $(\blb C,1)$.


\section{Exchange Representations}

Let $\lie g$ be a finite dimensional simple Lie algebra, and
$\hlie g$ the corresponding 
affine Lie algebra. Let $q\in {\mathbb C}^*,|q|\ne 1$, and let $U_q(\hlie g)$ be  
the associated quantum affine algebra 
(we will use the definition given in \cite{EFK98}). Let $\repf(U_q(\hlie g))$ 
be the category of finite dimensional (type 1) representations of $U_q(\hlie g)$.
For any $V \in \repf(U_q(\hlie g))$ and $z\in \blb C^*$ we can
consider the shifted representation $V(z)$ where all the Chevalley
generators act as before except for $e_0$ and $f_0$ which now act by
$ze_0$ and $z^{-1}f_0$ respectively.

We will need the following theorem. 

\begin{thm}\label{mero}  (Kazhdan-Soibelman, \cite{KaSo95}, Section 5)
Let $V,W\in\repf(U_q(\hlie g))$, and $\mathbf R_{VW}(z)$ be the series obtained by evaluating the universal R-matrix in $V(z)\otimes W$. Then $\mathbf R_{VW}(z)$ is an analytic function near zero, that extends meromorphically to $\blb C$. 
\end{thm}

We will give a new proof of this result (a bit simpler than the one in \cite{KaSo95}) in the appendix. 

Now let $\lie h$ denote a Cartan subalgebra of $\lie g$, and  
$\hlie h = \lie h \oplus \blb C$
the corresponding Cartan subalegbra in $\hlie g$. We denote the elements of $\hlie h^*$ by $\hat{\lambda} = (\lambda,k)$.
Then we can define the fusion matrix 
$$
J_{V,W}(u,\hat{\lambda}) : V(e^{-2\pi iu})\otimes W \to 
V(e^{-2\pi iu})\otimes W,
$$
which is a formal power series in $e^{2\pi iu}$
with coefficients being ${\rm End}(V\otimes W)$-valued 
rational functions of $q^{(\hat\lambda,\alpha_i)}$, where 
$\alpha_i$ are the simple roots. 
(see \cite{EtSc98b} and also \cite{EFK98}).
Namely, let $M_{\hat\lambda}$ be the Verma module over $U_q(\hlie
g)$ with highest weight $\hat\lambda$, and let $x_\lambda$ be its
highest weight vector. Let $v\in V$ be a homogeneous vector, and 
$\Phi_{\hat\lambda}^v(u):
M_{\hat\lambda}\to M_{\hat\lambda-{\rm weight}(v)}\hat\otimes V(z)$
be the intertwining operator such that the matrix coefficient 
$(x^*_{\hat\lambda-{\rm
    weight}(v)},\Phi_{\hat\lambda}^v(z)x_\lambda)$ equals $v$
(here $x_{\hat\mu}^*$ is the lowest weight vector of
$M_{\hat\mu}^*$). Then for any homogeneous vectors $v\in V$,
$w\in W$,
$$
J_{V,W}(u,\hat\lambda):=(x^*_{\hat\lambda-{\rm
    weight}(v)-{\rm weight}(w)},(\Phi^v_{\hat\lambda-{\rm
    weight}(w)}(e^{-2\pi iu})\otimes
1)\Phi^w_{\hat\lambda}(1)x_\lambda).
$$

Let ${\rm m}$ be the ratio of the squared lengths of long and short
roots of $\lie g$ (${\rm m}=1,2$ or $3$), and 
$h^\vee$ be the dual Coxeter number of $\lie g$.

\begin{prop}\label{p-adic}
The coefficients of the series 
$J_{V,W}(u,\hat\lambda)$ have an expansion into a power series 
in $p:=q^{-2{\rm m}(k+h^\vee)}$, 
and the sequence of these coefficients
tends to $0$ in the p-adic norm. 
\end{prop}

\begin{proof} This follows
from the ABRR equation for the fusion matrix \cite{EtSc99}
(which in our case is equivalent the qKZ equation of Frenkel and
Reshetikhin \cite{FrRe92},\cite{EFK98}).  
\end{proof}

Proposition \ref{p-adic} implies that one may consider 
the fusion operator $J_{V,W}(\hat\lambda):=J_{V,W}(0,\hat\lambda)$, which 
is a power series in $p$. 

\begin{prop}
The series $J_{V,W}(\hat\lambda)$ is convergent when 
$|p|<<1$, and the sum of this series 
extends to a meromorphic function in the region $\lbrace{(\lambda,k):
|p|<1\rbrace}$.  
\end{prop}

\begin{proof} This follows from the qKZ equation \cite{FrRe92},\cite{EFK98} and 
Theorem \ref{mero}.
\end{proof} 

Now define 
the exchange matrix with spectral parameter  
$$R_{V,W}(u,\hat{\lambda}) = J_{V,W}(u,\hat{\lambda})^{-1} 
\mathbf R_{W,V}^{21}(e^{2\pi i u}) J^{21}_{W,V}(-u,\hat{\lambda}).
$$
It follows by the arguments sketched above that 
both $J_{V,W}(u,\hat\lambda)$ and $R_{V,W}(u,\hat\lambda)$
are meromorphic functions on ${\mathbb C}\times \lie U$, 
where the open set 
$\lie U\subset \hlie h^*$ is defined by the condition 
$|q^{-2{\rm m}(k+h^\vee)}|<1$. Moreover, these arguments actually show that 
the divisor of singularities of these functions is 
a union of hyperplanes of the form $u-bk=c$, 
$<x,\lambda>-bk=c$, where $b,c\in \blb C, x\in \lie h$.

\begin{ex}\label{ex:sln}
Let $\lie g = \lie{sl}_n$. When $V = \blb C^n$, the formula for $R_{V,V}$ was obtained in \cite{Mou01}.  
\end{ex} 

\begin{thm}\label{thm:excisRm} \hfill
\begin{enumerate}\vspace{-.25cm} 
\item\cite{EtSc99} The fusion matrix satisfies the dynamical twist equation 
\begin{equation}\label{eq:dtwist}
J_{V\otimes W,U}^{12,3}(\hat\lambda) J_{V,W}^{12}(\hat\lambda - h^{(3)}) = J_{V,W\otimes U}^{1,23}(\hat\lambda)J_{W,U}^{23}(\hat\lambda)
\end{equation}

\item\cite{EtSc98b}  For any three finite dimensional representations 
$V_1,V_2,V_3$ of $U_q(\hlie g)$ the exchange matrix satisfies the
quantum dynamical Yang-Baxter equation 
\begin{align}\notag
R_{V_1,V_2}^{12}(u_1-u_2,\hat\lambda- h^{(3)}) R_{V_1,V_3}^{13}(u_1-u_3,\hat\lambda) 
R_{V_2,V_3}^{23}(u_2-u_3,\hat\lambda- h^{(1)})\\ \label{eq:qdybs1}
= R_{V_2,V_3}^{23}(u_2-u_3,\hat\lambda) R_{V_1,V_3}^{13}(u_1-u_3,\hat\lambda- h^{(2)}) 
R_{V_1,V_2}^{12}(u_1-u_2,\hat\lambda).
\end{align}
In particular, 
for any $V\in \repf(U_q(\hlie g))$ the function 
$R_{V,V}(u,\hat\lambda)$ is a quantum dynamical R-matrix with 
spectral parameter and step 1. 
\end{enumerate}
\end{thm}

\obs Strictly speaking, the last statement of the theorem contains an abuse 
of terminology, since $R_{V,V}$ is a meromorphic function not on $\blb C\times \hlie h^*$, 
but on its open subset $\blb C\times \lie U$. In this paper, we will ignore this 
terminological problem. 

Let $V$ be a finite dimensional representation of $U_q(\hlie g)$
and define $R_V = R_{V,V}$. It follows from theorem
\ref{thm:excisRm} that  for any $W \in \repf(U_q(\hlie g))$, 
the pair $(W,R_{V,W}(u,\hat{\lambda}))$ is a representation of $R_{V}$. It is called the {\bf exchange representation} of $R_V$ associated to $W$.

Now let us explain the main construction of this paper. 
Let $V \in \repf (U_q(\hlie g))$. Denote by $M$ the field of meromorphic functions on $\lie U$ and let $M_V \subset M$ be the subfield of functions
periodic under the weights of $V$. Set 
$\cal C(V) = \repf(U_q(\hlie g))\otimes_{\blb C} M_V$.
It is easy to see then that 
the assignment $W\to (W,R_{V,W})$ 
provides an exact, faithful functor $\cal F_V : \cal C(V) \to \repf(R_V)$.

One of our main results is the following theorem, 
which shows that for a nontrivial $V$, the category 
${\mathcal C}(V)$ (which is obtained from $\repf(U_q(\hlie g))$ 
by extension of scalars) is a full subcategory in $\repf(R_{V})$.  

\begin{thm}\label{main} The functor 
$\cal F_{V} : \cal C(V) \to \repf(R_{V})$ is a tensor functor.
Furthermore, if $V$ is nontrivial, then this functor is fully faithful. 
\end{thm}

\begin{proof}
We start with the following lemma. 

\begin{lem}\label{lem:tensoriso}
\begin{equation}
R_{V,W}^{12}(u,\hat\lambda - h^{(3)}) R_{V,U}^{13}(u,\hat\lambda) = J_{W,U}^{23}(\hat\lambda)^{-1}R_{V,W\otimes U}^{1,23}(u,\hat\lambda) J_{W,U}^{23}(\hat\lambda - h^{(1)})
\end{equation}
\end{lem}

\begin{proof} To shorten the notation we will write $V[u]$ instead of $V(e^{-2\pi iu})$. Let $\mathbf R_{W,V}^{21}[u] = \cal R^{21}|_{W\otimes V[u]}$ (analogously for other pairs of indices). Then
\begin{align*}
\text{l.h.s.} = & J_{V,W}^{12}(u,\hat\lambda - h^{(3)})^{-1}\mathbf R^{21}_{W,V}[u] J_{W,V}^{21}(-u,\hat\lambda-h^{(3)}) J_{V,U}^{13}(u,\hat\lambda)^{-1}\mathbf R^{31}_{U,V}[u] J^{31}_{U,V}(-u,\hat\lambda)  = \\
& J_{V[u],W}^{12}(\hat\lambda - h^{(3)})^{-1}\mathbf R^{21}_{W,V}[u] \Big(J_{W,V[u]}^{21}(\hat\lambda-h^{(3)}) J_{V[u],U}^{13}(\hat\lambda)^{-1}\Big)\mathbf R^{31}_{U,V}[u] J^{31}_{U,V[u]}(\hat\lambda)  \stackrel{\eqref{eq:dtwist}}{=} \\ 
& J_{V[u],W}^{12}(\hat\lambda - h^{(3)})^{-1}\mathbf R^{21}_{W,V}[u] \Big(J_{W\otimes V[u],U}^{21,3}(\hat\lambda)^{-1} J_{W,V[u]\otimes U}^{2,13}(\hat\lambda)\Big)\mathbf R^{31}_{U,V}[u] J^{31}_{U,V[u]}(\hat\lambda) =\\  
& \Big(J_{V[u],W}^{12}(\hat\lambda - h^{(3)})^{-1} J_{V[u]\otimes W,U}^{12,3}(\hat\lambda)^{-1}\Big)\mathbf R^{21}_{W,V}[u]\mathbf R^{31}_{U,V}[u]\Big( J_{W,U\otimes V[u]}^{2,31}(\hat\lambda)  J^{31}_{U,V[u]}(\hat\lambda)\Big) \stackrel{\eqref{eq:dtwist}}{=} \\
& \Big(J_{W,U}^{23}(\hat\lambda)^{-1} J_{V[u],W\otimes U}^{1,23}(\hat\lambda)^{-1}\Big) \mathbf R^{21}_{W,V}[u]\mathbf R^{31}_{U,V}[u]\Big( J_{W\otimes U,V[u]}^{23,1}(\hat\lambda)  J^{23}_{W,U}(\hat\lambda-h^{(1)})\Big) = \text{r.h.s.} 
\end{align*}
The fourth equality is clear while the last follows from the hexagon axiom for the R-matrix.
\end{proof}

\obs A more general version of this lemma in the non-affine situation is used in 
\cite{EtVa00}.

Now let us prove that $\cal F$ is a tensor functor.  
To do this, define the tensor structure \linebreak $\cal F_V(W)\odot \cal F_V(U) \to \cal F_V(W\otimes U)$ to be given by
$$J_{W,U}(\hat{\lambda}) : (W,R_{V,W})\odot (U,R_{V,U}) \to (W\otimes U,R_{V,W\otimes U})$$
It follows from Lemma \ref{lem:tensoriso} that $J_{W,U}$ is a
morphism (in fact an isomorphism) and from \eqref{eq:dtwist} that
it satisfies the properties of a tensor structure on a functor,
as required. 

The rest of the section is dedicated to the proof of the full faithfulness
of $\cal F$ in the case when $V$ is nontrivial.

Following the conventions of \cite{EFK98}, set $L_W^+(z):=\mathbf
R_{W,V}^{21}(z)$. 

\begin{lem}\label{lem:M(W,U)}
One has ${\rm Hom}_{U_q(\hlie g)}(W,U) = \{ A\in \text{Hom}_{U_q(\hlie h)}(W,U)$; $L^{+}_U(z)^{-1}(1\otimes A)L^{+}_W (z)$ is of zero weight in both components$\}$. 
\end{lem}

\begin{proof}
It is clear that the first space is a subspace 
in the second one. For the converse it is enough to prove that, if $A:W \to U$ is 
such that $L^{+}_U(z)^{-1}(1\otimes A)L^{+}_W (z)$ is of zero
weight in both components, then $A \in \text{Hom}_{U_q(\hlie
  g)}(W,U)$. To do this, recall that 
the universal R-matrix on modules with central charge $0$ can be written as
$$\cal R = q^{\sum x_i\otimes x_i}\big(1 + \sum_{i\geq 0} (q_i - q_i^{-1}) e_i\otimes f_i +\dots \big)$$
where $x_i$ form an orthonormal basis for $\lie h$ and the
remaining terms are in $U_q(\hlie n_+)\otimes U_q(\hlie n_-)$ and
correspond to other weights. 
This means that for nontrivial $V$, the condition that the
expression 
$L^{+}_U(z)^{-1}(1\otimes
A)L^{+}_W (z)$ is of zero weight in both components implies
that $A$ is of zero weight, and 
$f_i|_V\otimes[e_i,A] = 0 \foral i\ge 0$ (one just needs to look at the terms corresponding to simple roots in the expression of the universal R-matrix). It follows that $[e_i,A]=0$. Hence $A$ is an interwiner over $U_q(\hlie b_+)$.
But it is proved in Proposition 4.1.3 of \cite{KaSo95} 
that the images of $U_q(\hlie b_+)$ and $U_q(\hlie g)$ 
in the endomorphism algebra of any finite dimensional
$U_q(\hlie g)$-module are the same. Hence, $A$ is an intertwiner 
 over the whole $U_q(\hlie g)$, as desired.  
\end{proof}

Let $\cal M(W,U):={\rm Hom}_{U_q(\hlie g)}(W,U)\otimes_{\blb C} M_V$
be the space of morphisms between $W$ and $U$ in the category $\cal C(V)$, and let 
$\widetilde{\cal M}(W,U):={\cal M}(W,U)\otimes_{M_V}M$. Lemma \ref{lem:M(W,U)} 
implies that $\widetilde{\cal M}(W,U)$ is the space of solutions over $M$
of the system of linear algebraic equations in finitely many variables
(entries of $A$), expressing the condition that 
$L^{+}_U(z)^{-1}(1\otimes A)L^{+}_W (z)$ is of zero weight in both components. 

Now let ${\cal N}(W,U)={\rm Hom}_{\repf(R_V)}({\cal F}(W),{\cal F}(U))$.
Below we will show that the space 
$\widetilde{\cal N}(W,U):={\cal N}(W,U)\otimes_{M_V}M$ is contained in 
the space of solutions over $M$ of a {\bf deformation} of this system 
of linear algebraic equations. This will imply that 
${\rm dim}{\cal N}(W,U)\le {\rm dim}{\cal M}(W,U)$, 
which immediately implies that the injective map 
${\cal M}(W,U)\to {\cal N}(W,U)$ defined by ${\cal F}$ is actually an isomorphism, 
i.e., ${\cal F}$ is fully faithful. 

To do this, we need the following lemma. 

\begin{lem}\label{lem:M'(W,U)}
Let $L_W, L_U$ be the exchange operators $R_{V,W}, R_{V,U}$ and 
$\widetilde{\cal M}\,'(W,U) = \{ A(\hat\lambda) \in \text{Hom}_{U_q(\hlie h)}(W,U)\otimes_{\blb C}M; L_U(u,\hat\lambda)^{-1}(1\otimes A(\hat\lambda))L_W (u,\hat\lambda)$ is of zero weight in both components$\}$. Then 
the natural map $\xi:{\cal N}(W,U)\otimes_{M_V}M \to \widetilde{\cal M}\,'(W,U)$
is injective.
\end{lem} 

\begin{proof}
We need to show that any $M_V$-linearly independent set $\{a_1, \dots, a_m\} \subset {\cal N}(W,U)$ is also linearly independent over $M$. Assume it is false
and take a counterexample with the smallest length. We can assume that $a_m = g_1 a_1 + \dots + g_{m-1} a_{m-1}$ with $g_i \in M$ but not all of them in $M_V$. Then
\begin{align*}
a_m(\hat\lambda - h^{(1)}) & = L_U(u,\hat\lambda)^{-1}(1\otimes a_m(\hat\lambda))L_W (u,\hat\lambda)\\ & = \sum g_i(\hat\lambda) L_U(u,\hat\lambda)^{-1}(1\otimes a_i(\hat\lambda))L_W (u,\hat\lambda) = \sum g_i(\hat\lambda) a_i(\hat\lambda -h^{(1)})
\end{align*}
Let $v \in V$ be homogeneous of weight $\mu$. Then, applying the last equality to 
$v\otimes {\rm Id}$, we have
$$a_m(\hat\lambda-\mu) = \sum g_i(\hat\lambda)a_i(\hat\lambda -\mu) \quad\Leftrightarrow\quad a_m(\hat\lambda) = \sum g_i(\hat\lambda+\mu)a_i(\hat\lambda)$$
Hence, for some $\mu$, $\sum \big(g_i(\hat\lambda)-g_i(\hat\lambda+\mu)\big)a_i(\hat\lambda)$ is a non-trivial (since not all $g_i$ are in $M_V$) 
vanishing linear combination of an even smaller length. Contradiction.
\end{proof}

Now it remains to prove that $L_W(u,\hat\lambda)$ is a deformation of 
$L_W^{+}(z)$, where $z=e^{2\pi iu}$.

Let $\kappa = k+h^\vee{}$ and $p=q^{-2{\rm m}\kappa}$. We will write $k \to \infty$ to mean $p \to 0$ keeping $q$ fixed. 

\begin{lem}\label{lem:klimit}
Let  $J_{V,W}(\lambda)$ be the fusion matrix corresponding to the the ``finite dimensional'' quantum group $U_q(\lie g)$. Then $\lim_{k\to \infty} J_{V,W}(u,\hat\lambda) = 
J_{V,W}(\lambda)$.
\end{lem}

\begin{proof}
Let $v$ and $w$ be homogeneous vectors in $V$ and $W$
respectively and consider the correlation function \cite{EFK98}
$\psi(z,\hat\lambda) = z^{\Delta}J_{V,W}(u,\hat\lambda)(v\otimes
w)$ where $\Delta$ depends on the weights of $v$ and $w$ and on
$\hat\lambda$. 
It is easy to see that 
$$\psi(z,\hat\lambda) = z^{\Delta}\sum_{m\geq 0} \psi_m(\hat\lambda)z^m$$
and that $\psi_0(\hat\lambda) = J_{V,W}(\lambda)(v\otimes w)$. Furthermore, $\psi$ satisfies the q-\kz  equations, from which we get
$$\psi(pz,\hat\lambda) = \mathbf R_{W,V}^{21}(pz) q_{(2)}^{\Lambda} \psi(z,\hat\lambda)$$
where $\Lambda$ depends only on $\lambda$ and the weights of $v,w$ (not on $k$). Writing $\mathbf R_{W,V}^{21}(z) = \sum_{j\geq 0} \mathbf R^{21}_j z^j$ we find
$$p^{\Delta} \sum_{m\geq 0}\psi_m(\hat\lambda) p^m z^m = \Big(\sum_{j\geq 0} \mathbf R^{21}_j p^j z^j\Big) q^{\Lambda}_{(2)}\Big(\sum_{i\geq 0} \psi_i(\hat\lambda)z^i\Big)$$
Actually $p^{\Delta}$ is some constant $c$ depending only on $q,\lambda$ and the weights of $v,w$. We collect the term in $z^l$ in both sides to get
$$(cp^l - \mathbf R_0^{21}q^{\Lambda}_{(2)}) \psi_l = \sum_{m=1}^l \mathbf R^{21}_m p^m q^{\Lambda}_{(2)} \psi_{l-m}$$
from where we conclude by induction that $\psi_l \to 0$ as $k \to \infty$ for all $l>1$, as stated.
\end{proof}

Now assume that $|q|<1$ and let $\lambda \to -\infty$ mean ${\rm
  Re}(\lambda,\alpha_i) \to -\infty$, where $\alpha_i$ are the simple
roots. Then, by \cite{EtVa00}, $J_{V,W}(\lambda) \to 1$ as $\lambda \to
-\infty$, and we conclude that $L_W(u,\hat\lambda) \to
\mathbf R_{W,V}^{21}(z)$ 
and, therefore, $L_W(u,\hat\lambda)$ is a deformation of
$L_W^{+}(z)$.
One can proceed analogously for $|q|>1$. Theorem \ref{main} is proved. 
\end{proof}


\section{Gauge Transformations}

To prove the theorem stated in the introduction, it remains to 
identify the category $\repf(R_V)$ for $V$ being the vector
representation of $U_q({\hlie sl}_n)$ with the category 
of representations of Felder's elliptic quantum group. 
To construct such an identification is the goal of this section. 
It is achieved by using the fact that the corresponding 
dynamical R-matrices are gauge equivalent \cite{Mou01}.

The concept of gauge transformations of quantum dynamical R-matrices of $\lie{gl}_n$ type, that we briefly recall now, was introduced in \cite{EtVa98}. 

Let $\lie h$ be an abelian $n$-dimensional Lie algebra. A multiplicative $m$-form on $\lie h^*$ is a collection of meromorphic functions $\varphi = \{\varphi_{a_1,\dots,a_m}(\lambda)\}$ where $\{a_1, \dots, a_m\}$ runs through all (ordered) subsets of $\{1,\dots,n\}$ such that
$$\varphi_{a_1,\dots,a_{i+1},a_i,\dots,a_m}(\lambda)\varphi_{a_1,\dots,a_m}(\lambda)=1$$
Let $\Omega^m$ be the set of all multiplicative $m$-forms. With the obvious definitions one turns $\Omega^m$ into an abelian group where the neutral element is the constant form $\mathbf{1}_{a_1,\dots,a_m}(\lambda) = 1$. 

For a given $\gamma \in \blb C^*$ and $s \in \{1, \dots, n\}$ define the operator $\delta_s$ on the space of meromorphic functions on $\lie h^*$ given by
$$\delta_s f(\lambda_1, \dots, \lambda_n) = \frac{f(\lambda_1, \dots, \lambda_n)}{f(\lambda_1, \dots,\lambda_s-\gamma,\dots, \lambda_n)}$$
and then the differential homomorphism $d_{\gamma} : \Omega^m \to \Omega^{m+1}$, $\varphi \mapsto d_{\gamma}\varphi$, where
$$(d_{\gamma}\varphi)_{a_1,\dots, a_{m+1}}(\lambda) = \prod_{s=1}^{m+1} \big(\delta_{a_s} \varphi_{a_1,\dots,a_{s-1},a_{s+1},\dots, a_{m+1}}(\lambda)\big)^{(-1)^{s+1}}$$
We have $d^2_{\gamma}\varphi =\mathbf{1} \foral\varphi$. A form $\varphi$ is said to be $\gamma$-closed if $d_{\gamma}\varphi = \mathbf{1}$ and $\gamma$-exact if $\varphi = d_{\gamma}\psi$ for some form $\psi$. When $\gamma =1$ we just say that $\varphi$ is closed or exact and denote $d_1$ simply by $d$.

We list now the gauge transformations that will be important for us. Let $R$ be an R-matrix of $\lie{gl}_n$ type and step $\gamma$ and define the following transformations

\begin{equation}\label{eq:g4}
R(u,\lambda) \mapsto c(u)R(u,\lambda)
\end{equation}
where $c(u)$ is a meromorphic scalar function;
\begin{equation}\label{eq:g5}
R(u,\lambda) \mapsto R(au,b\lambda + \mu)
\end{equation}
where $a,b \in \blb C^*$ and $\mu \in \lie h^*$;
\begin{gather}\notag
R(u,\lambda) \mapsto \\\label{eq:g2}
\sum_{m=0}^{n} \alpha_{m,m}(u,\lambda)E_{m,m}\otimes E_{m,m} + \sum_{m\neq l} \varphi_{m,l}(\lambda)\alpha_{m,l}(u,\lambda) E_{m,m}\otimes E_{l,l} + \beta_{m,l}(u,\lambda) E_{l,m}\otimes E_{m,l}
\end{gather}
where $\varphi=\{\varphi_{m,l}(\lambda)\}$ is a $\gamma$-closed 2-form.

It was proved in \cite{EtVa98} that if $\tilde R$ is obtained from $R$ by these gauge transformations, then it is a quantum dynamical R-matrix of $\lie{gl}_n$ type and step $\gamma$ (for \eqref{eq:g4} and \eqref{eq:g2}) or $\gamma/b$ (for \eqref{eq:g5}). It is clear that if $\tilde R$ is obtained from $R$ by \eqref{eq:g4} and \eqref{eq:g5}, then $\repf(\tilde R)$ is equivalent to $\repf(R)$. The following theorem, concerning gauge transformation \eqref{eq:g2}, was proved in \cite{EtVa99} in the language of bialgebroids. We repeat the proof here with the appropriate adaptations for this context.

\begin{thm}\label{thm:g2}
Let $R$ be a quantum dynamical R-matrix of $\lie{gl}_n$ type and step $\gamma$ and $\varphi$ be a $\gamma$-exact 2-form. If $\tilde R$ is obtained from $R$ by \eqref{eq:g2} with this $\varphi$ then $\repf(\tilde R) \cong \repf(R)$
\end{thm}

\begin{proof}
Let $\varphi = d_{\gamma}\zeta$ and set $\xi = \sum_a \zeta_a E_{a,a}$. Note first that we can write
$$\tilde R(u,\lambda) = (\xi^{(1)}(\lambda-\gamma h^{(2)}))^{-1} (\xi^{(2)}(\lambda))^{-1} R(u,\lambda) \xi^{(1)}(\lambda) \xi^{(2)}(\lambda-\gamma h^{(1)})$$
Then let $(W,L)$ be a representation of $R$ and set
$$\tilde L(u,\lambda) = (\xi^{(1)}(\lambda-\gamma h^{(2)}))^{-1}L(u,\lambda)\xi^{(1)}(\lambda)$$
Now it is easy to check that the functor $(W,L) \mapsto (W,\tilde L)$ (identity on morphisms) is an equivalence of tensor categories.
\end{proof}

One can regard the exchange matrix $R_{\blb C^n}$ of example \ref{ex:sln} as defined on the whole Cartan sub-algebra of $\lie{gl}_n$. Then, if one looks at the coordinate $k$ as a parameter, it is a quantum dynamical R-matrix of $\lie{gl}_n$ type.  

\begin{thm}\label{thm:Rk=Rell}\cite{Mou01}
Let $R^{ell}_{\tau,\gamma}$ be the elliptic quantum dynamical R-matrix of example \ref{ex:felder}. Then $R_{\blb C^n}$ is gauge equivalent to $R^{ell}_{\tau,\gamma}$ for $\gamma = \frac{\log q}{\pi i}$ and $\tau = -\kappa\gamma = -\kappa\frac{\log q}{\pi i}$, where $\kappa = k+n$.
\end{thm}

The three gauge transformations we have listed are used in theorem \ref{thm:Rk=Rell}. In particular, the gauge transformation of type \eqref{eq:g2} used is given by the following multiplicative 2-form 
\begin{equation}\label{eq:varphi}
\varphi_{i,j}(\lambda) = \sigma_{j,i}(\lambda/\gamma-\rho)
\end{equation}
where, for $m<l$,
$$\sigma_{l,m}(\lambda) = q\frac{\Gamma_p(1+\frac{1}{\kappa}(\lambda+\rho)_{l,m}+\frac{1}{\kappa})} {\Gamma_p(1+\frac{1}{\kappa}(\lambda+\rho)_{l,m})} \frac{\Gamma_p(-\frac{1}{\kappa}(\lambda+\rho)_{l,m})} {\Gamma_p(-\frac{1}{\kappa}(\lambda+\rho)_{l,m}+\frac{1}{\kappa})} \quad, \sigma_{m,l}(\lambda) = \frac{1}{\sigma_{l,m}(\lambda)},$$
and $\Gamma_q$ is the q-Gamma function. 

\begin{lem}\label{lem:exact}
The multiplicative 2-form \eqref{eq:varphi} is $\gamma$-exact.
\end{lem}

\begin{proof}
It is enough to show that $\tilde{\varphi}_{i,j}(\lambda) = \varphi_{i,j}(\gamma\lambda)$ is exact. Let $\xi = \{\xi_j\}$ be the 1-form given by $\xi_j(\lambda) = \prod_{i<j}q^{\lambda_i}$. Then $\delta_s\xi_j(\lambda) = 1$ if $s>j$ and $\delta_s\xi_j(\lambda) = q$ if $s<j$. Therefore, for $m<l$, $(d\xi)_{m,l}(\lambda) =q$ and $(d\xi)_{l,m}(\lambda) = q^{-1}$. In a similar way one proves that $\tilde{\varphi} = d(\xi\eta\zeta)$ where 
$$\eta_j(\lambda) = \prod_{i<j}\Gamma_p\Big(\frac{\lambda_{i,j}+1}{\kappa}\Big)^{-1} \qquad\text{and}\qquad \zeta_j(\lambda) = \prod_{i<j}\Gamma_p\Big(1+\frac{\lambda_{j,i}}{\kappa}\Big)^{-1}$$ 
\end{proof}

Thus, we have constructed the equivalence of categories ${\cal G}: \repf(R_V)\to \repf(R_{\tau,\gamma}^{ell})$, where $V=\blb C^n$, and hence, proved the following theorem, which is our main
result. 

\begin{thm}\label{mainth} The functor ${\cal F}={\cal G}
\cal F_V:\repf(U_q(\hlie{sl}_n))\otimes_{\blb C}M_V\to \repf(R^{ell}_{\tau,\gamma})$
is an exact, fully faithful tensor functor, which defines an
equivalence of $\repf(U_q({\hlie sl}_n))\otimes_{\blb C}M_V$ with a full subcategory of
$\repf(R^{ell}_{\tau,\gamma})$.
\end{thm}


\section{Evaluation Representations}

In this section we would like to study the action of 
the functor $\cal F$ on evaluation representations of $U_q(\hlie
sl_n)$. Let $W_\nu(z)$ be the evaluation representation of 
$U_q({\hlie sl}_n)$ with dominant integral highest weight $\nu$, 
evaluated at $z\in \blb C^*$. 

\begin{prop} The representation ${\cal F}(W_\nu(z))$ 
is isomorphic to the evaluation
representation $V_\mu(z)$ of the elliptic quantum group defined 
in \cite{TaVa01} (see Corollary 3.4 and Theorem 5.9),
tensored with a 1-dimensional representation. 
\end{prop}

\begin{proof} Similarly to quantum affine algebra,
for elliptic quantum group there is a natural notion of a highest weight
for a finite dimensional representation $W$ (it was studied in detail
in \cite{Ca01}). Namely, the highest weight is, essentially, the set
of eigenvalues of the diagonal entries of the 
noncommutative matrix $L_W(u,\hat\lambda)$ on the highest weight
vector. 

It is not difficult to compute highest weights of 
representations ${\cal F}(W_\nu(z))$ and show that they 
coincide with highest weights of $V_\nu(z)$ up to tensoring with
a 1-dimensional representation. Namely one needs to prove this 
statement for the fundamental representations (exterior powers),
since for other representations the statement follows by
taking tensor products. 

To prove the statement for exterior powers, one first proves 
it for the vector representation $V$, in which case it is trivial,
since it reduces to analyzing entries of $R_V$. To pass to an 
arbitrary exterior power, one may use fusion procedure (i.e.
the fact that $\Lambda^iV$ is contained in the product 
of $i$ copies of $V$ with appropriate shifts). 

Now, it was shown in \cite{Ca01} that if two irreducible 
highest weight finite dimensional representations of the 
elliptic quantum group
have the same highest weight, then they are isomorphic. 
This completes the proof. 
\end{proof}

\obs It is clear that the functor $\mathcal F$ extends 
naturally to a tensor functor $\text{\bf Rep}_{\cal O}(U_q({\hlie
  sl}_n))\otimes_{\blb C}M_V\to 
\text{\bf Rep}_{\cal O}(R_{\tau,\gamma}^{ell})$ between the
categories $\mathcal O$ of representations of the quantum affine algebra,
respectively the elliptic quantum group
(i.e. categories of representations whose $\lie h$-weight 
multiplicities are finite, and the multiplicity function 
is supported on a finite union of sets of the form 
$\mu-Q_+$, where $Q_+$ is the positive cone in the root lattice
of $\lie g$). Namely, the construction of $\mathcal F$ we gave
generalizes to this case in a straighforward way.
Moreover, it is not hard to see that after this extension the
functor remains fully faithful (the proof of full faithfulness 
also generalizes without significant changes). 
In particular, we see that any evaluation module for 
$U_q({\hlie sl}_n)$ admits an elliptic deformation. 
This was conjectured in \cite{TaVa01}, and the conjecture was
checked for Verma modules and finite dimensional representations.


\section{Comparison of $\mathcal F$ with the classical
  Kazhdan-Lusztig functor}

The sense in which the functor $\cal F$ is a
q-analog of the Kazhdan-Lusztig functor is not
straightforward, and needs to be explained. 

Indeed, the Kazhdan-Lusztig functor is 
a functor from the fusion category of representations of an
affine Kac-Moody algebra $\hlie g$ at some level $k$ to the category 
of representations of the quantum group $U_v(\lie g)$, where 
$v=e^{\frac{\pi i}{{\rm m}(k+h^\vee)}}$. On the other hand,
 in this paper we do not seem to be  
considering any q-fusion categories. 

Nevertheless, there 
is a sense in which $\cal F$ is a q-analogue of the
Kazhdan-Lusztig functor. A two-line justification is that while the Kazhdan-Lusztig 
functor arises from the analysis of the monodromy of the
Knizhnik-Zamolodchikov equations, the functor $\cal F$ arises
in a similar way from the analysis of the monodromy of the quantum
Knizhnik-Zamolodchikov equations. To some extent this is explained 
in \cite{FrRe92} and \cite{KaSo95}. However, let us give   
a somewhat more vivid explanation, in the 
special case $\lie g=\frak{sl}_n$.
This explanation uses the small elliptic quantum group of Tarasov
and Varchenko. 

Recall that the classical Kazhdan-Lusztig functor $\cal F_0$, essentially,
does the following. 
Given a finite dimensional representation $V$ 
of $\lie g$ and a (generic) level $k$, it 
introduces on $V$ an action of the quantum group $U_v(\lie g)$, 
with $v=e^{\frac{\pi i}{{\rm m}(k+h^\vee)}}$. More precisely, 
one first uses the exchange construction for $\hlie g$ at level
$k$ to turn $V$ into a representation of the exchange dynamical
quantum group $F_v(G)$, $G=SL_n$ (see \cite{EtVa99}), and then turns it into a representation
of $U_v(\lie g)$ by sending the dynamical parameter to infinity. 

On the other hand, given a (generic) level $k$
and a finite dimensional representation $V$ of $U_q(\lie g)$,
 the functor $\cal F$ introduces on $V$ an action 
of the elliptic quantum group 
$E_{\tau,\gamma/2}(\lie g)$, with $\tau=-\frac{(k+h^\vee)\log
  q}{\pi i}$, and $\gamma=\frac{\log q}{\pi i}$ (by first affinizing $V$ and then
applying the functor $\cal F$).  

Now, one may check using modular
invariance of theta-functions that $E_{\tau,\gamma/2}(\lie g)$ is 
isomorphic to $E_{-1/\tau,\gamma/2\tau}(\lie g)$ (the
corresponding dynamical R-matrices are gauge equivalent). 
Thus, $V$ becomes a representation 
of $E_{-1/\tau,\gamma/2\tau}(\lie g)$. This representation is easily seen to be
an evaluation representation, i.e., it factors through the small 
elliptic quantum group $e_{-1/\tau,\gamma/2\tau}(\lie g)$ (\cite{TaVa01}).

Finally, we claim that 
if $k$ is fixed then the algebra 
$e_{-1/\tau,\gamma/2\tau}(\lie g)$ 
is a q-deformation of the exchange dynamical quantum group $F_{v}(G)$
(which allows us to regard $\cal F$ as a q-deformation of $\cal F_0$). 
Indeed, it is easy to check that as $\omega\to \infty$, 
$e_{\omega,a}(\lie g)$ degenerates into 
$F_{e^{-2\pi ia}}(G)$.  
Since we have $-1/\tau=\frac{\pi i}{(k+h^\vee)\log q}$, and 
$\gamma/2\tau=-\frac{1}{2(k+h^\vee)}$, the claim follows. 

\appendix
\section{Appendix : Meromorphicity of the R-matrix}

In this appendix we give a new proof of Theorem \ref{mero}.

Let $U$ be a neighborhood of zero in $\blb C^n$, and $G: U\to \blb C^n$ be an analytic mapping, such that $G(0)=0$. Let $p\in \blb C, |p|>1$.

\begin{lem}\label{estimate}
Let $f\in z\blb C^n[[z]]$ be a formal solution to the difference equation 
$$ f(pz)=G(f(z)) $$
Then $f$ converges to an analytic function in some neighborhood of zero.
\end{lem}

\begin{proof}
Let $f(z)=\sum_{n\ge 1}f_n z^n$, and 
$G(y)=\sum_{n\ge 1}g_ny^n$, where $g_n$ are multilinear operators. Since $G$ is analytic at zero, there exists $A>0$ such that the norms of the coefficients are estimated by $|g_n|<A^n$ for all $n$. 

Let $k_0$ be so big that for $k\ge k_0$ one has$|(p^k-g_1)^{-1}|<2|p|^{-k}$, and $2A<|p|^{k/2}$. Let $C$ be such that $|f_k|<C$ for $k<k_0$, and let $B>\text{max}(1,\frac{AC}{|p|^{1/2}-1})$. 

Let us prove by induction in $k$ that $|f_k|<CB^{k-1}$. The statement is clear for $k<k_0$ (as $B>1$). Assume that the statement is known for all $n<k$, for some $k\ge k_0$, and prove it for $k$. We have
$$ f_k=(p^k-g_1)^{-1}\sum_{1<r\le k}\sum_{i_1,...,i_r:\sum i_j=k} g_r(f_{i_1},...,f_{i_r}) $$
Now, $|g_r(f_{i_1},...,f_{i_r})|
\le |g_r|\prod_{l=1}^r |f_{i_l}|<(AC)^rB^{k-r}$, and the number
of terms with $g_r$ is dominated by the binomial coefficient
$\left(\begin{matrix} k-1\\ r-1\end{matrix}\right)$ (as $i_j\ge 1$, $r>1$). Thus, we have  
$$ |f_k|<2(AC/B)(1+AC/B)^{k-1}|p|^{-k}B^k< CB^{k-1}, $$
since $2A<|p|^{k/2}$ and $1+AC/B<|p|^{1/2}$. 
The inductive step is proved. 

This implies that $f(z)$ converges for small $z$. We are done.
\end{proof}

Now we proceed with the proof of Theorem \ref{mero}.
Recall the crossing symmetry condition (equation (9.19) of \cite{EFK98})  
\begin{equation}\label{eq:cs}
(((\mathbf R_{V,W}(z)^{-1})^{t_1})^{-1})^{t_1}) =
(q^{2\rho}\otimes 1)\mathbf R_{V,W}(zq^{2{\rm m}h^\vee{}})(q^{-2\rho}\otimes 1)
\end{equation}
where $t_1$ is transposition in the first component and $\rho$ is
the half-sum of positive roots 
of $\lie g$. It is valid for any finite dimensional
representations. If $|q|\ne 1$, it can be regarded as an equation
of the type $f(pz)=G(f(z))$ considered above.
Namely, $f$ takes values in the
space of endomorphisms of $V\otimes W$,
$p=q^{\pm 2{\rm m}h^\vee}$, 
$G(X)=(((X^{-1})^{t_1})^{-1})^{t_1}$, 
and the origin in $\text{End}(V\otimes W)$ is at the point $X=\mathbf R_{V,W}(0)$, 
which is a fixed point of $G$. Thus, 
by Lemma \ref{estimate}, 
the function $\mathbf R_{V,W}(z)$ is analytic near zero. Then the meromorphicity in the whole complex plane follows immediately from the difference equation \eqref{eq:cs} and the fact that in our case 
$G$ is a birational isomorphism.\hfill$\square$ 

\obs It is clear that the same result and proof applies to $V,W$ being in category \cal O as $U_q(g)$ modules (since the R-matrix has zero weight, and the weight spaces of $V\otimes W$ are finite dimensional).



\end{document}